\documentclass{IEEEtran4PSCC}
\usepackage{nomencl}
\makenomenclature
\usepackage{etoolbox}
\usepackage[pdftex,hypertexnames=false,linktocpage=true]{hyperref}
\hypersetup{colorlinks=true,linkcolor=blue,anchorcolor=blue,citecolor=blue,filecolor=blue,urlcolor=blue,bookmarksnumbered=true,pdfview=FitB}
\renewcommand\nomgroup[1]{%
  \item[\bfseries
  \ifstrequal{#1}{A}{Indices and Graph}{%
  \ifstrequal{#1}{B}{Model Parameters}{%
  \ifstrequal{#1}{C}{Decision Variables}{}}}%
]}
\newtheorem{remark}{Remark}
\usepackage{paralist}
\usepackage[ruled,vlined]{algorithm2e}
\ifCLASSINFOpdf
   \usepackage[pdftex]{graphicx}
\else
   \usepackage[dvips]{graphicx}
\fi
\usepackage[cmex10]{amsmath}
\usepackage[dvipsnames]{xcolor}
\usepackage{threeparttable}

\usepackage{soul}

\makeatletter
\let\old@ps@headings\ps@headings
\let\old@ps@IEEEtitlepagestyle\ps@IEEEtitlepagestyle
\def\psccfooter#1{%
    \def\ps@headings{%
        \old@ps@headings%
        \def\@oddfoot{\strut\hfill#1\hfill\strut}%
        \def\@evenfoot{\strut\hfill#1\hfill\strut}%
    }%
    \def\ps@IEEEtitlepagestyle{%
        \old@ps@IEEEtitlepagestyle%
        \def\@oddfoot{\strut\hfill#1\hfill\strut}%
        \def\@evenfoot{\strut\hfill#1\hfill\strut}%
    }%
    \ps@headings%
}
\makeatother

\psccfooter{%
        \parbox{\textwidth}{\hrulefill \\ \small{23nd Power Systems Computation Conference} \hfill \begin{minipage}{0.2\textwidth}\centering \vspace*{4pt} \includegraphics[scale=0.06]{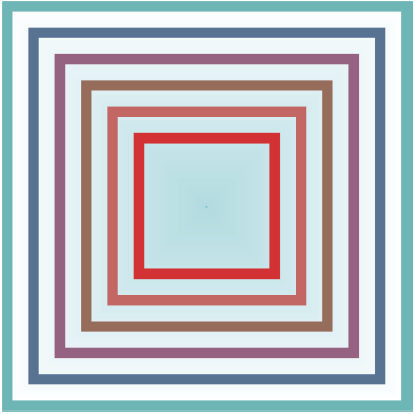}\\\small{PSCC 2024} \end{minipage} \hfill \small{Paris-Saclay, France --- June 4 -- June 7, 2024}}%
}

\begin{document}

\setstcolor{purple}

%
\title{An OpenStreetMaps based tool to study the energy demand and emissions impact of electrification of medium and heavy-duty freight trucks}

\author{
\IEEEauthorblockN{Nawaf Nazir\\ Bowen Huang\\ Shant M. Mahserejian}
\IEEEauthorblockA{Pacific Northwest National Laboratory, Richland, USA \\
\{nawaf.nazir, bowen.h, shant.mahserejian\}@pnnl.gov}
}


\maketitle

\begin{abstract}

In this paper, we present the mathematical formulation of an OpenStreetMaps (OSM) based tool that compares the costs and emissions of long-haul medium and heavy-duty (M\&HD) electric and diesel freight trucks, and determines the spatial distribution of added energy demand due to M\&HD EVs. The optimization utilizes a combination of information on routes from OSM, utility rate design data across the United States, and freight volume data, to determine these values. In order to deal with the computational complexity of this problem, we formulate the problem as a convex optimization problem that is scalable to a large geographic area. In our analysis, we further evaluate various scenarios of utility rate design (energy charges) and EV penetration rate across different geographic regions and their impact on the operating cost and emissions of the freight trucks. Our approach determines the net emissions reduction benefits of freight electrification by considering the primary energy source in different regions. Such analysis will provide insights to policy makers in designing utility rates for electric vehicle supply equipment (EVSE) operators depending upon the specific geographic region and to electric utilities in deciding infrastructure upgrades based on the spatial distribution of the added energy demand of M\&HD EVs. To showcase the results, a case study for the U.S. state of Texas is conducted.

\end{abstract}


\thanksto{\noindent The authors acknowledge support from the Laboratory Directed Research Development (LDRD) program at PNNL.}

  \nomenclature[A]{$\mathbf{G}$}{ full set of locations in the graph}
  \nomenclature[A]{$\mathbf{A}$}{ all possible undirected arcs connecting origin and destination of $\mathbf{G}$}
  \nomenclature[A]{$k$}{ index for vertices/locations in the graph}
  \nomenclature[A]{$t$}{ index for time period, $t = 1, \cdots, T$}

  \nomenclature[B]{$d_{k}$}{ travel distance on arc $(k,k+1), \forall k \in\mathbf{G}$}
  \nomenclature[B]{$t_{k}$}{ travel time on arc $(k,k+1), \forall k \in\mathbf{G}$}
  \nomenclature[B]{$v_{k}$}{ travel speed on arc $(k,k+1), \forall k \in\mathbf{G}$}
  \nomenclature[B]{$\lambda_{k}$}{ utility energy rate price at location $k$ (\$/MWh)}
  \nomenclature[B]{$\eta^+$}{ battery charging efficiency}
  \nomenclature[B]{$\eta^-$}{ battery discharging efficiency}
  \nomenclature[B]{$\eta_\text{w2e}$}{ average tractive energy power at wheels by battery and motor (kW/miles)}
  \nomenclature[B]{$\overline E$}{ maximum battery capacity (kWh) of vehicle available to charge}
  \nomenclature[B]{$\underline E$}{ minimum battery capacity (kWh) of vehicle available for discharging}
  \nomenclature[B]{$\overline P$}{ maximum charging power for vehicle (kW)}
  \nomenclature[B]{$\underline P$}{ minimum charging power for vehicle (kW)}
  \nomenclature[B]{$e_0$}{ initial SOC of the vehicle battery at the beginning of the day/time horizon}
  \nomenclature[B]{$e_{K+1}$}{ terminal SOC of the vehicle battery at the end of the day/time horizon}
  \nomenclature[B]{$Q_m$}{ maximum capacity (tons) of vehicle}

  \nomenclature[C]{$x_{k}$}{ binary variable representing charging decision at location $k, \forall k\in\mathbf{G}$, if vehicle charge at location $k, x_{k}=1$, otherwise $x_{k}=0$}
  \nomenclature[C]{$e_k^+$}{ energy charged at location $k$}
  \nomenclature[C]{$e_{(k,k+1)}^-$}{ energy consumption on arc $(k,k+1)$}
  \nomenclature[C]{$p_k^+$}{ battery charging power at location $k$}
  \nomenclature[C]{$p_k^-$}{ battery discharging power on arc $(k,k+1)$}
  \nomenclature[C]{$l_{k}$}{ vehicle load (tons) on arc $(k,k+1)$}
  \nomenclature[C]{$r$}{ the potential route between given origin and destination pair (O-D pair)}
\printnomenclature

\section{Introduction}
\IEEEPARstart{I}n recent years, the electrification of vehicle fleets has emerged as a promising solution to mitigate environmental concerns and reduce reliance on fossil fuels in transportation, a sector which accounts for nearly 30\%~\cite{FFGHG} of all GHG emissions. As societies worldwide strive to transition to sustainable and greener modes of transportation, fleet electrification, particularly for long-haul freight trucks, has gained significant attention from both industry and academia. This shift towards electric vehicles (EVs) in fleet operations offers numerous potential benefits, including reduced carbon emissions, lower operational costs, and enhanced energy efficiency. However, the electrification of long-haul freight trucks, despite its potential benefits, presents unique challenges. These vehicles are integral to the global logistics and supply chain network, often covering extensive distances and requiring high power demands. Consequently, electrifying this specific type of fleet is confronted with issues related to battery capacity, charging infrastructure, and the need for cost-effective solutions to ensure long-haul capabilities while minimizing environmental impact. A critical aspect in evaluating the economic viability of fleet electrification lies in understanding the impacts of utility rate design on EV charging patterns and figuring out the spatial distribution of the resulting charging energy demand. Stakeholders ranging from utility companies, which can plan infrastructure upgrades based on this distribution, to policy makers who can design effective regulations, stand to benefit significantly from this insightful analysis. Efficient location based charging is pivotal for fleet operations as it directly influences transportation costs, energy consumption, and overall fleet performance. Therefore, conducting an economic assessment of fleet electrification through vehicle charging optimization becomes imperative to accurately evaluate the financial feasibility, emissions reduction, spatial distribution of energy demand and potential advantages associated with transitioning to EVs.

Numerous studies have focused on the economic assessment of fleet electrification~\cite{tamba2022economy,kleiner2015electrification,martinez2021assessment}, specifically incorporating vehicle routing optimization techniques~\cite{yang2015electric,lin2016electric,xiao2021electric,ham2021electric}. These studies, spanning a range of scales from urban area to statewide, have explored the integration of optimization models and algorithms with fleet electrification considerations to facilitate decision-making processes and provide valuable insights for fleet operators and policymakers.

One common approach in the literature involves developing mathematical models that incorporate various factors related to fleet electrification and vehicle routing. These models often consider parameters such as vehicle range, battery capacity, charging infrastructure availability~\cite{worley2012simultaneous}, energy consumption~\cite{sassi2014vehicle}, and charging time~\cite{yang2015electric}. By combining these elements into optimization frameworks, the existing works aim to determine the optimal allocation of EVs across delivery routes, charging strategies, and scheduling decisions to minimize costs and maximize efficiency. 

Moreover, researchers have also investigated the impact of fleet electrification on total cost of ownership (TCO)~\cite{kumar2020total}, taking into account factors such as vehicle acquisition costs, maintenance expenses, energy costs, and potential incentives or subsidies~\cite{hagman2016total}. These TCO models, combined with vehicle routing optimization in ~\cite{jahic2021impact}, enable a comprehensive assessment of the economic viability of fleet electrification, providing insights into the potential cost savings and payback periods associated with the transition to EVs.

Furthermore, advancements in data analytics and real-time data collection in ~\cite{sevilla2022state} have allowed for the integration of dynamic factors into fleet electrification economic assessments. By incorporating real-time traffic data, weather conditions, and demand fluctuations, \cite{song2012adaptive} aims to develop more accurate and adaptive optimization models that account for real-world uncertainties and dynamic operational environments.

However, it's noteworthy that certain critical aspects have been underexplored in the literature. Notably, very few studies addressing long-haul medium and heavy-duty freight truck fleets, where electrification potential intersects with significant logistical and energy demand complexities. Additionally, comprehensive emissions reduction calculations are often lacking, along with an understanding of the spatial distribution of the resulting charging energy demand, accounting for utility rate design. Moreover, the need for scalable solutions that can be applied to larger geographic regions has not received adequate attention. The electrification of such fleets demands formulating a highly complex, large-scale optimization problem that may not be computationally tractable to solve, given the complexity of interconnected variables and constraints. To bridge these gaps, we present our solution, which extends the existing body of research by addressing these challenges in a holistic manner, providing valuable insights and scalable optimization techniques to support the electrification of long-haul medium and heavy-duty freight truck fleets.

In this paper, we aim to provide valuable insights and an applicable tools for fleet operators, policymakers, and stakeholders in understanding the financial implications, benefits, and challenges associated with transitioning to electric vehicles for long-haul medium and heavy duty freight trucks. 
What sets our research apart is our comprehensive approach to addressing the complexities of electrifying these critical freight fleets. We combine utility rate design, freight volume data, and route data from OpenStreetMaps (OSM) into a comprehensive optimization framework. Unlike previous studies, we formulate this optimization problem in a scalable manner, allowing for practical application over extensive geographic regions, from state wide to nationwide. This scalability enables us to address challenges presented by the electrification of long-haul freight truck fleets, where multiple variables and constraints interact within a vast operational landscape.

The structure of this paper will be as follows, in section~\ref{sec:formulation}, we describe the problem in the mathematical formulation, provide detailed constraints and assumptions and then propose a numerical algorithm framework to process the user input and integrate the mathematical formulation into the optimization engine, hence solve a large scale vehicle routing problem iteratively. In section~\ref{sec:data}, the data flow management and preprocessing will be shown along with the user interface. In section~\ref{sec:results}, we explore the scalability of the algorithm by solving state-wide transportation system in Texas, compare the environmental and economic metrics and access the efficiency of our framework. Finally the conclusion are laid out in Sec.~\ref{sec:conclusion} .

\section{Problem formulation}\label{sec:formulation}
In this section, we outline the research methodology. Following the problem description, we describe the primary objective of our problem. Subsequently, we will remark and detail our assumptions on the key elements of our optimization framework, including the incorporation of utility rate design, freight volume data, and route data from OSM.

\subsection{System Description}
Consider a transportation system represented by the graph $\mathbf G$, comprised of fleet charging at designated locations, denoted by $k=1,2,\ldots, K\in \mathbf{G}$, the set of undirected arcs connecting the origin $O\in\mathbf{G}$ and destination $D\in\mathbf{G}$ (O-D pair) is denoted by $\mathbf A$. To find the optimal route, charging location and charging power of the given vehicle fleet, the input data for the optimization problem includes the energy price at location $k$, $\lambda_k$, all the routes between the given O-D pair and the energy charged at location $k$, $e^+_{k,r}$.
\subsection{Objective}
In the proposed EV Scheduling Problem, given pairs of origin and destination and the graph involving all possible arcs connecting them, we are minimizing the charging cost of the EV fleet while satisfying the daily duty. The optimal decision includes the charging location (where to charge), charging power (how fast to charge) and the optimal objective value will determine the optimal route (which route to select).
\begin{align}
    \text{minimize} \quad \sum_{r=1}^R\sum_{k=1}^{K}\lambda_{k} e_{k,r}^+\label{eq:obj1}
\end{align}
The process is repeated to obtain aggregate values over a span of one year. The rationale behind choosing a yearly timeframe is multifold. Firstly, the Bureau of Transportation Statistics (BTS) datasets~\cite{FAF5} is released on an annual basis, providing a comprehensive view of transportation metrics over a 12-month period. This allows for consistent and reliable comparisons, using a standardized dataset. Additionally, considering an entire year ensures that variances across different months and seasons are captured. Seasonal differences can significantly affect transportation patterns due to factors like weather conditions, holidays, and industrial demands. By encompassing all these variances, a yearly analysis offers a more holistic and representative insight into transportation trends and impacts.

\subsection{Constraints}
\subsubsection{Battery capacity constraints}
\begin{eqnarray}
    \underline E \leq e_0 + \sum_{k=1}^K e_k^+\eta^+ - \sum_{k=1}^K \frac{e_{(k,k+1)}^-}{\eta^-} \leq \overline E,&\; \label{batcst1}\\
    x_k\left[\underline E + \max{(\frac{e_{k-1}^-}{\eta^-},\frac{e_{k}^-}{\eta^-})}\right]\leq e_k^+\eta^+\leq x_k\overline E,&\label{batcst2}
\end{eqnarray}
where $ \forall k,k-1 \in\mathbf{G}$, the eq.~\eqref{batcst1} ensures that energy consumption on battery for the entire trip between origin and destination(undirected) is bounded by energy capacity limits. For the charging energy between location $k$ and $k+1$, the eq.~\eqref{batcst2} enforces the battery charging energy to be safe in the worst cases, i.e., charging energy to be able to arrive at nearest location. $x_k$ is the binary variable to decide whether or not to charge at the location $k$.\\
\subsubsection{Load constraints}
\begin{eqnarray}
    e_{(k,k+1)}^- = p_k^- t_k \quad \forall k\in\mathbf{G}\label{loadcst1}\\
    p_k^- = \eta_\text{w2e}v_k\frac{l_k}{Q_m}, \forall k\in \mathbf{G}\label{loadcst2}
\end{eqnarray}
The eq.~\eqref{loadcst1} ensures that EV  on the arc $(k,k+1)$ discharged enough energy to meet the load and travel demand. The eq.~\eqref{loadcst2} utilize the tractive energy efficiency from the wheels to the battery to model the discharge power during the trip as function of load level and speed.\\


\subsubsection{Other constraints}
\begin{eqnarray}
    d_{k} = v_{k} t_{k},\quad \forall k\in\mathbf{G}\label{travelcst}\\
    x_{k}\in\{0,\;1\},\quad \forall k \in \mathbf{G}\label{binarycst}
\end{eqnarray}
The eq.~\eqref{travelcst} shows the relationship between travel distance, travel time and speed. Eq.~\eqref{binarycst} defines the binary decision variable $x_{k}$. \\
Hence the overall optimization problem is formulated as,
\begin{align}
    \text{minimize} \quad& \sum_{r=1}^R\sum_{k=1}^{K}\lambda_{k} e_{k,r}^+\label{eq:obj1}\\
    \text{subject to}\quad&\text{eq.}\eqref{batcst1}-\eqref{binarycst}
\end{align}

Note that by introducing the binary variable for the battery charging and discharging state, the optimization problem becomes a mixed-integer linear programming (MILP) problem. Tackling large-scale MILP problems, however, is known to be computationally challenging due to their inherent complexity, especially as the problem size grows~\cite{zafar2023decentralized}. Our proposed approach aims not only to model the problem but also to present a algorithm framework that is scalable, addressing the computational challenges typically associated with large-scale MILP problems.

To facilitate the integration process of the proposed framework with our input data, we also made assumptions on the EV fleet charging and discharging to relax some restrictions for the optimization problem.
\begin{remark}
In our approach, we adopt a 24-hour look-ahead window to optimize the charging locations consecutively throughout the entire year and for the entire fleet. This approach allows for dynamic, day-to-day adaptability, capturing the fluctuating demands of fleet electrification. Instead of tackling the computational challenge of optimizing the whole year in one iteration, this daily sequential method divides the yearly problem into smaller sub-problems. This not only enhances computational efficiency but also provides a granular perspective that aligns with real-world dynamics.
\end{remark}
\begin{remark}
Consider a single pair of origin and destination points at a time. This method greatly simplifies our optimization task. Our fleet's characteristics – long distances, fixed start and end points, set loads, daily routines, specific charging parameters, and SOC transitions – allow us to view potential routes as undirected. This aligns with the BTS freight volume data, which details origin-destination pairs, enabling parallel optimization for each pair. Iterating over these routes ensures scalability thorough our analysis, optimizing our algorithm's performance within the BTS data framework.
\end{remark}

To further ensure the convexity of our optimization problem formulation, we supplement with the following assumptions:
\begin{itemize}
    \item Travel speed on each arc/edge is constant but may vary across arcs.
    \item Battery charging time is fixed and does not affect the 24 hour driving schedule.
    \item All the vertices including origin and destination are available for charging.
    \item All the vehicles associated with given origin and destination pair have the same schedule of charging and route selection.
    \item The total load demand served on a route does not exceed the vehicle fleet capacity.
\end{itemize}

\subsection{Algorithm framework}\label{subsec:algorithm}
\begin{figure}[!htp]
\centering
\includegraphics[width=0.44\textwidth]{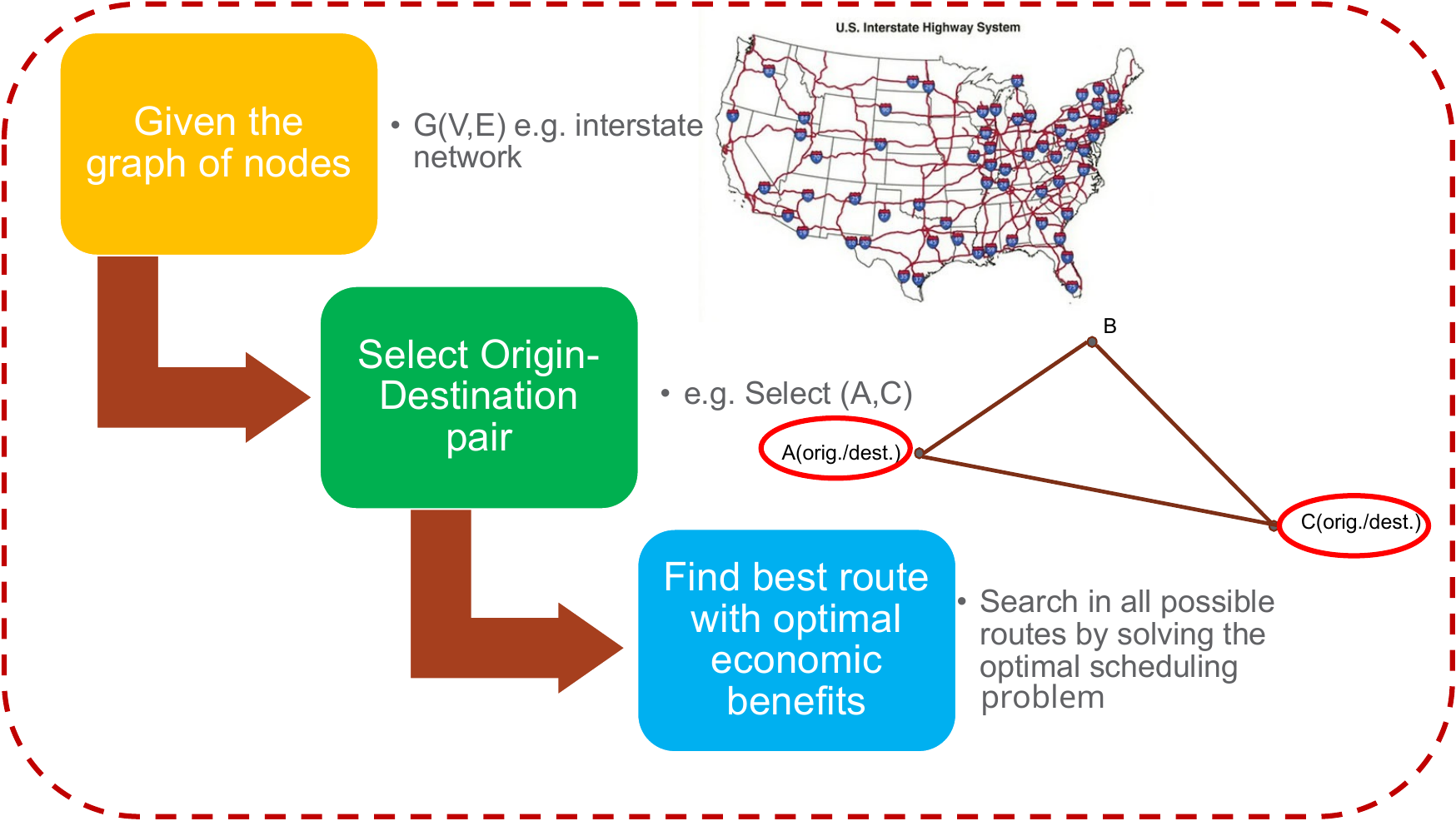}
\caption{EV routing and dispatch algorithm framework used to calculate the optimal charging location and charging power and determine the spatial distribution of added energy demand.}
\label{fig:algorithm_blk}
\end{figure}

Here we elaborate the framework of the algorithm to solve the EV optimal dispatch and routing problem defined above.
The input parameter of the proposed algorithm consists of the following:
\begin{itemize}
    \item Fleet Information: The algorithm requires input data on the fleet size, types of vehicles (including EVs), their charging capacities, and vehicle operating costs.
    \item Charging Infrastructure: The availability, locations, and charging capacities of charging stations need to be considered to determine the feasibility of EV routing and charging plans.
    \item Electricity Pricing: The algorithm requires information on electricity pricing, including peak and off-peak rates, to optimize the charging schedules of EVs.
\end{itemize}
Given the data input from either the fleet owner or the utility company to optimize their infrastructures, the algorithm framework will take use of the \textit{\bf Julia} engine as the optimizer and go through the steps iteratively in the Fig.~\ref{fig:algorithm_blk}.
The detailed algorithm flow is shown in Algorithm~\ref{alg:efect}.
\begin{algorithm}[h]
\SetAlgoLined
 \LinesNumbered
 Consider the list of routes between given undirected O-D pair, $r=1,2,\ldots,R$, $t=1,2,\ldots,365$,\\
    \For{\textnormal{route }$r=1,2,\ldots,R$}{
       \For{\textnormal{day }$t=1,2,\ldots,365$}{
       Solve the refined optimization problem for the selected route and O-D pair,\label{alg:s1}
    \begin{align*}
        \text{minimize} \quad &\sum_{k=1}^{K}\lambda_{k} e_{k}^+\\
        \text{subject to}\quad&\text{eq.}\eqref{batcst1}-\eqref{binarycst}
    \end{align*}}
    Calculate the optimal energy cost, $\text{CO}_2$ emissions, and energy demand for the entire fleet.\label{alg:s2}}
    Find the route with the most cost saving and emission reduction, return the optimal route and charging location solution.
    \caption{Optimal routing and charging location algorithm for the fleet electrification}\label{alg:efect}
\end{algorithm}

\section{Data Management and Preparation}\label{sec:data}
In order to ensure practical and meaningful results for end-use stakeholders of the developed tool, the following two real-world open data source were consolidated along a real-world representation of possible fleet routes: flow of freight volume along routes, and local energy pricing along routes. The route development took place prior to any analysis, which had two dependencies itself, namely the cities that can serve as origins and destinations appropriate for fleet vehicle routes, and the road network map along which they can travel. To interface with a potential user stakeholder, a user interface (UI), called the EFECT-UI, was developed to collect route and fleet information, and then to execute corresponding optimization calculations. 

\begin{figure}[!t]
\centering
\includegraphics[width=0.42\textwidth]{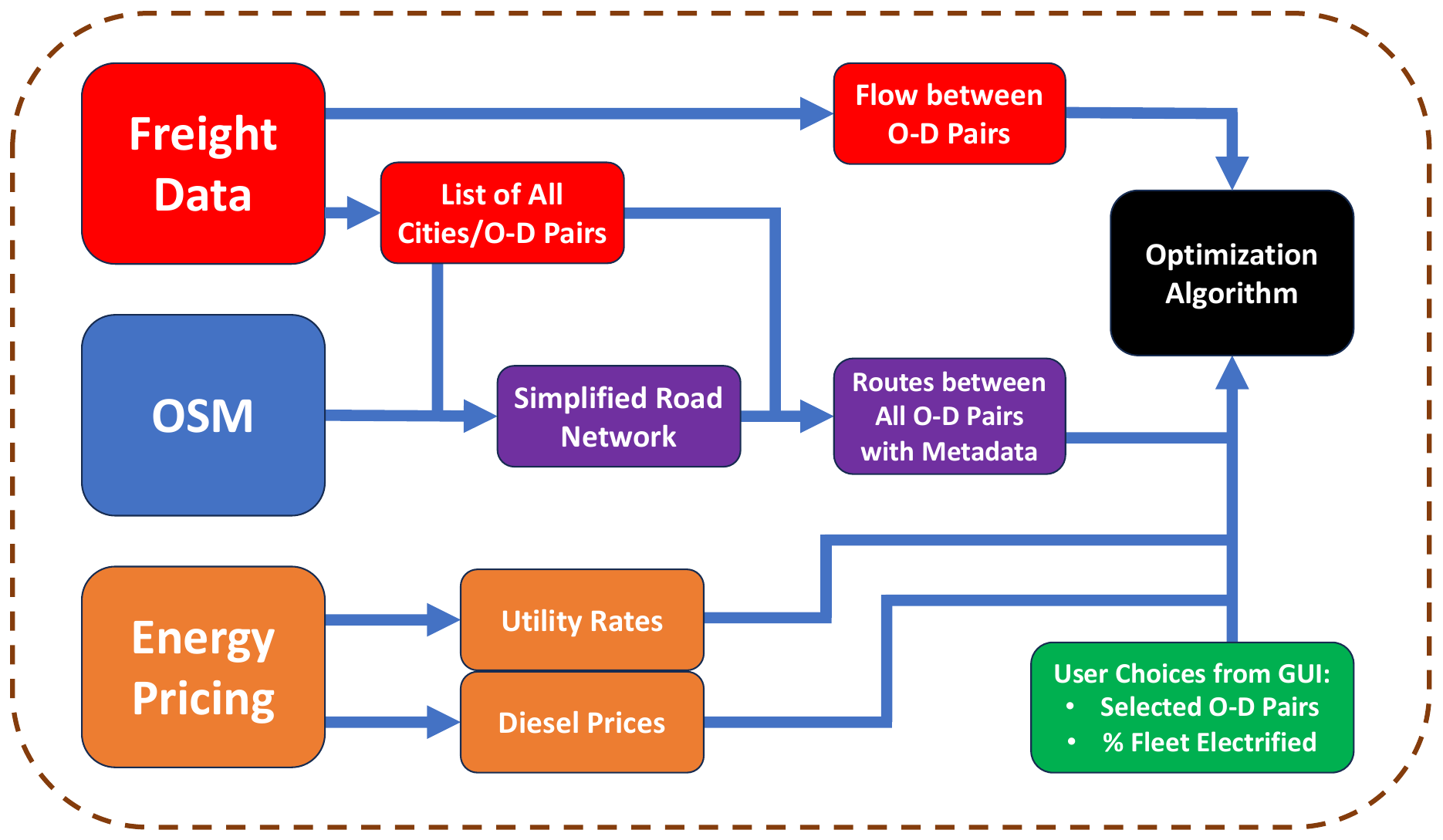}
\caption{Flow of information from data sources, and how they are prepared for input into the optimization algorithm.}
\label{fig:dataflow}
\end{figure}

\begin{figure}[!t]
\centering
\includegraphics[width=0.44\textwidth]{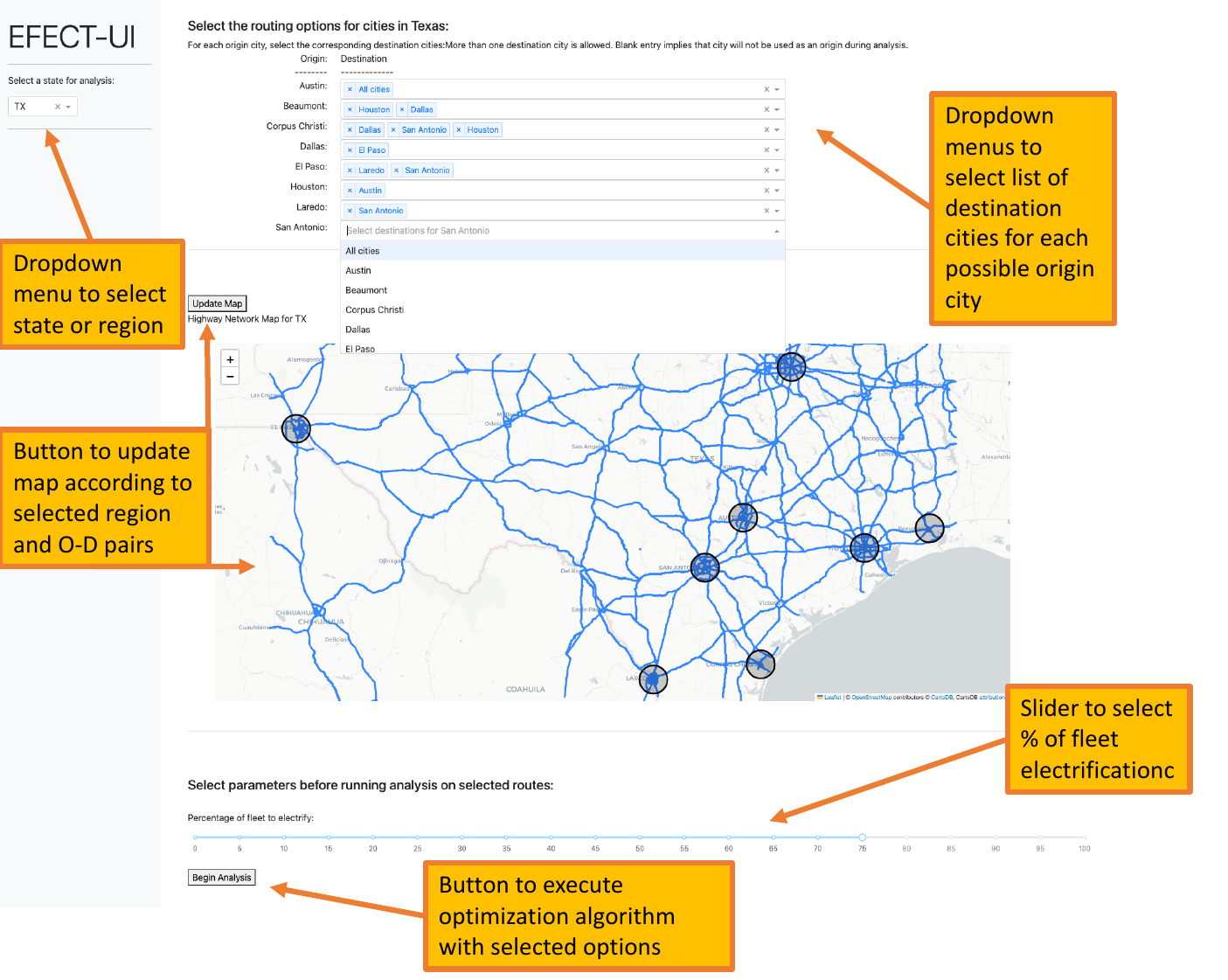}
\caption{A screenshot of the EFECT-UI developed to identify optimal charging locations, estimates of total energy costs, and potential reduction to emissions. Users select a region, cities for O-D pairs, and analysis parameters before executing the algorithm. Maps for the road network and selected routes are displayed.}
\label{fig:gui}
\end{figure}

This section covers the details of the flow of information, displayed in Figure \ref{fig:dataflow}, where data from original sources are processed to develop routes and their corresponding metadata, and they are combined with user-selected routes and parameters to provide input into the optimization algorithm. Figure \ref{fig:gui} displays a screenshot of the EFECT-UI representing the developed tool, where a user can select the geographical area bounded by a US state or other relevant region, cities in that region to serve as origin-destination (O-D) pairs for fleet routes, and additional parameters, such as percentage of a fleet to be electrified. 

To develop routes, we first relied on data supplied by the Freight Analysis Framework Version 5 (FAF5) prepared by US Department of Transportation (DOT) Federal Highway Administration (FWHA) BTS \cite{FAF5}. This data captures comprehensive freight movement between states and major metropolitan areas by integrating data from multiple sources. We made use of information corresponding to the "Truck" mode of transportation for this study. The included list of cities were used to limit O-D pairs to serve as starting and ending points for routes. The FAF5 dataset included flow (in tons) between each O-D pair, which we also used to model the traffic between different city pairs. The FAF5 dataset is represented as the red boxes in Figure \ref{fig:dataflow}.

To build the road network, we relied on Open Street Maps (OSM), an open-source distribution of geographic data for the world, including detailed road connections and location based metadata \cite{OpenStreetMap}. Latitude-longitude coordinates and corresponding counties of network node locations served particularly useful for our study when associating potential charging sites to the correct local utilities responsible for supplying the required energy. We used a combination of the \texttt{osmnx} and \texttt{networkx} packages in Python to first generate a localized map of detailed streets around each metro city area, and then combined the major highway network only throughout a buffered bounding box defined by the group of cities included for a given region.

Once the road network map was defined, the fastest routes between each O-D pair was generated, and potential en-route charging sites were identified at approximately 50~km increments. The routes were organized into JSON files, which each route was named using an indexed O-D pairs and were supplemented with metadata related to distance, time duration, human interpretable highway names, and charging site location information. The road network development is represented with the purple boxes in Figure \ref{fig:dataflow}.

The energy pricing data consists of the utility rate data and diesel price data. For the battery electric vehicle (BEV) fleet, we consider the utility rate data from the The Utility Rate Database (URDB) hosted by OPENEI~\cite{openei}. For the internal combustion engine vehicle (ICEV) fleet, the diesel price data at state level is obtained from American Automobile Association (AAA)'s real time data for 2023~\cite{aaaprice}. This information is represented by the orange boxes in Figure \ref{fig:dataflow}

The EFECT-UI integrates the proposed algorithm framework and the data management functionality into the user interface (UI), a snapshot of EFECT-UI is shown in Fig.~\ref{fig:gui}. The interface also includes some visualization of the road network and selected routes, which makes use of the list of network edges for each route, was also included in the JSON files created. Based on user selections, a subset of this JSON file serves as the input into the algorithm described in in section~\ref{subsec:algorithm}. Contributions from user-chosen criteria and values are displayed as the green box and connected lines in Figure \ref{fig:dataflow}, as they are aggregated with the rest of the data before being fed into the optimization algorithm.

\section{Case studies and analysis}\label{sec:results}
\begin{figure}[!t]
\centering
\includegraphics[width=0.44\textwidth]{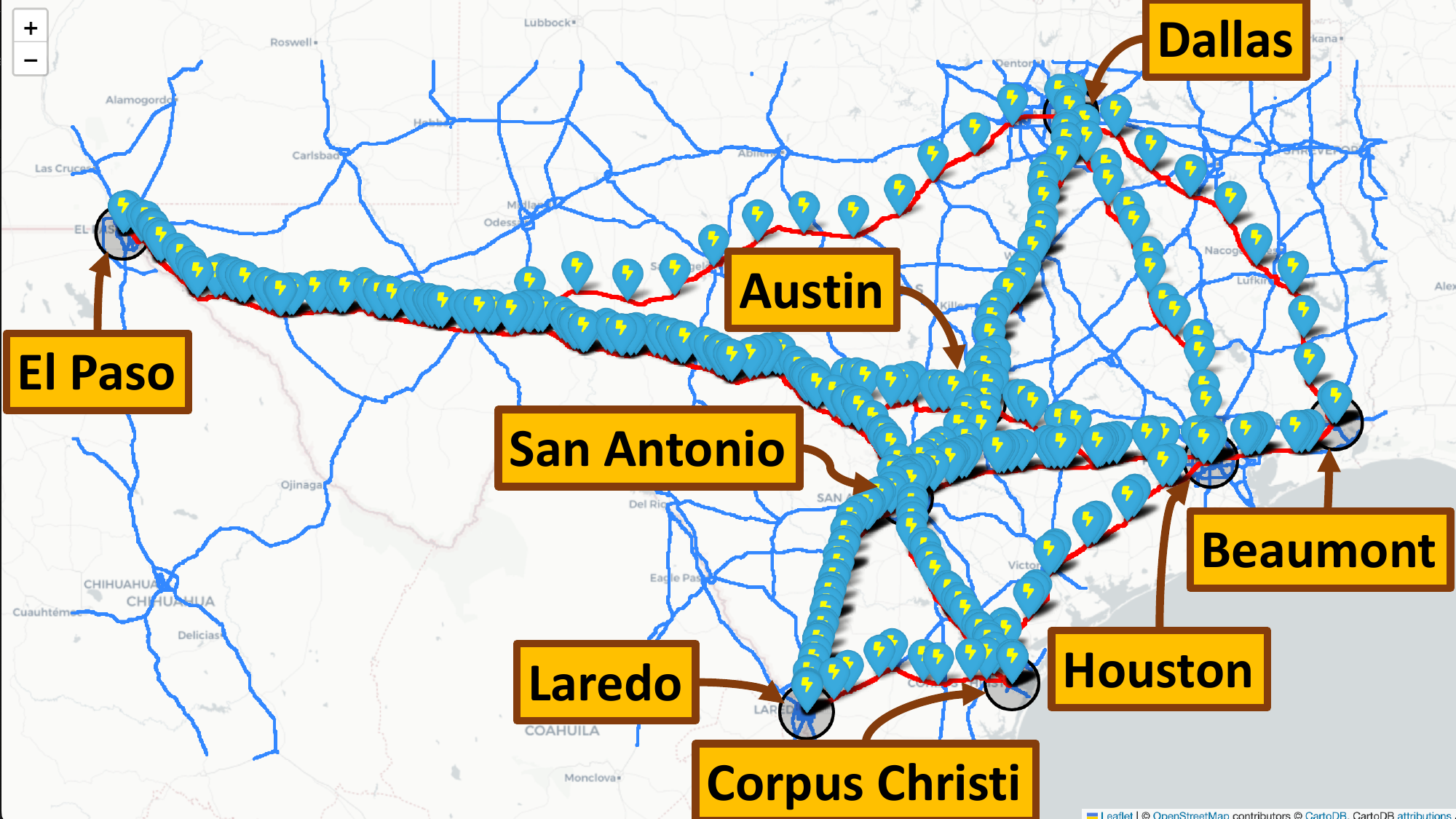}
\caption{Texas road network map (blue), with routes between all possible origin-destination pairs (red) and possible charging sites (markers).}
\label{fig:texas_routes}
\end{figure}

To showcase the results and applicability of our algorithm framework, we selected the state of Texas as a representative case study. Within Texas, key cities such as Austin, Beaumont, Corpus Christi, Dallas, El Paso, Houston, Laredo, and San Antonio have been considered. The detailed road network, illustrating the routes and connections between these cities, is presented in Figure \ref{fig:texas_routes}.

The same techniques can be scaled up to include regions defined by regional transmission organization (RTO) boundaries or at the national level. It's worth mentioning that the Texas state boundaries are analogous to the Electric Reliability Council of Texas (ERCOT) RTO, and though technically El Paso is not included within ERCOT boundaries, the potential charging locations en-route are, and thus including El Paso among the possible O-D pairs list is necessary for a scope around ERCOT. 

To simplify the process of solving the problem, we first generate the list of potential fastest routes between the given O-D pairs from the OSM data then solve the optimal charging solution in the optimization engine for the entire year.

The case study is considering the 8 cities in Texas, and selects all the possible O-D pairs among these 8 cities to approximate the total energy cost (in dollar values), total emission (in tons of $\text{CO}_2$ equivalent) and total energy consumptions. 
Emission calculations for both ICEV and BEV fleets stem from their primary energy sources
:
ICEVs: Emissions originate from fuel combustion. They are calculated by multiplication of amount of diesel consumption into diesel factor.
BEVs: Emissions hinge on the electricity grid's composition. They are calculated by multiplication of annual carbon intensity from \cite{energysourcemap}. The primary energy sources are taken into account for the carbon emission intensity from different locations in Texas.
The detailed cost and emission comparison between BEV and ICEV fleet over entire state is shown in \autoref{tab:comparison}.
\begin{table}[!t]
    \centering
    \begin{threeparttable}
      \caption{Total Operational Cost and Emission Comparison between BEV and ICEV fleets for the state of Texas.} \label{tab:comparison} \centering
      \renewcommand{\arraystretch}{1.1}
      \setlength{\tabcolsep}{2.5pt}
        \begin{tabular}{cccccc}
          \hline
          Type of Fleet  & Cost(\$)       & Emission(100 tons $\text{CO}_2$ eq) & Energy(MWh)\tnote{2} \\ \hline
          BEV            & 2.3M\tnote{1}  & 202                    & 10975            \\ \hline
          ICEV           & 25.6M          & 2141                   & 116324         \\ \hline
        \end{tabular}
    \begin{tablenotes}
        \small{
        \item[1] M: million
        \item[2] 1 gallon U.S. diesel oil = 0.0407 MWh}
    \end{tablenotes}
    \end{threeparttable}
\end{table}
Based on the optimal charging schedule for the entire year and the entire fleet, given the freight volume load between the O-D pairs, the BEV fleet's energy cost (operational cost only) can  be reduced to about 9\% of the ICEV fleet's cost. While ICEV fleet emits over 10 times the $\text{CO}_2$ equivalent emissions of BEVs, underscoring the significant environmental benefits of electric fleets, especially in terms of greenhouse gas emissions. The energy consumption for the ICEV fleet is approximately 10.6 times that of the BEV fleet. This substantial difference underscores the inefficiencies inherent in ICEVs, given that a significant portion of the energy in fuel is lost as heat rather than being used for propulsion.

\begin{figure}[!htp]
\centering
\includegraphics[width=0.5\textwidth]{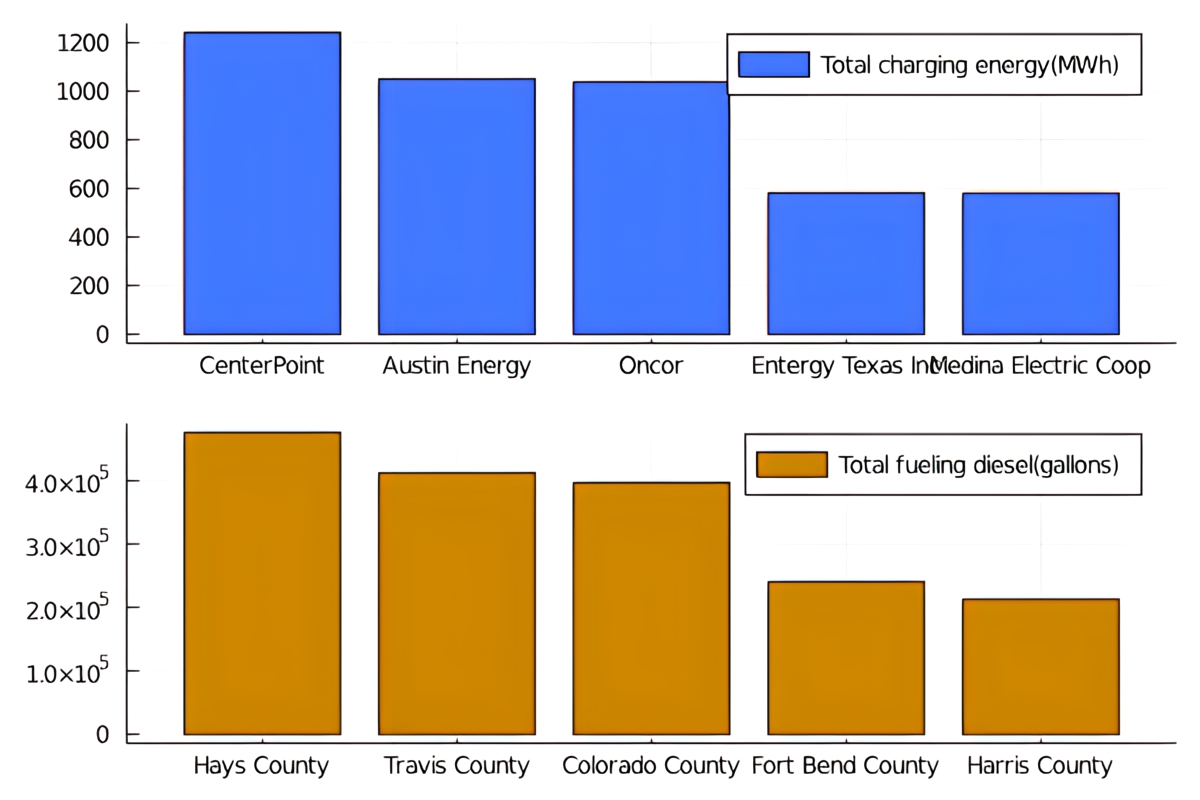}
\caption{Energy demand from BEV fleet at utility level and Fuel demand from ICEV fleet at county level for the state of Texas.}
\label{fig:demand_bar}
\end{figure}

The spatial energy demand distribution comparison between the BEV fleet and ICEV fleet, is shown in Fig.~\ref{fig:demand_bar}. The charging demand for utility companies along x axis is shown on the top for BEV fleet, and the fueling demand for ICEV fleet is shown in the bottom at county level. For the top half of Fig.~\ref{fig:demand_bar}, it illustrates the charging demand of BEVs for various utility companies. The CenterPoint, Austin Energy and Oncor provides the top 3 energy demand locations over the selected routes. For the bottom half of figure, we list the amount of diesel fuels consumed by the ICEV fleet, since we assign the fueling location and diesel price only to the county level. In the figure we can observe that Hays County provides the most fueling demand followed by Travis County and Colorado County.

\begin{figure}[!htp]
\centering
\includegraphics[width=0.75\linewidth]{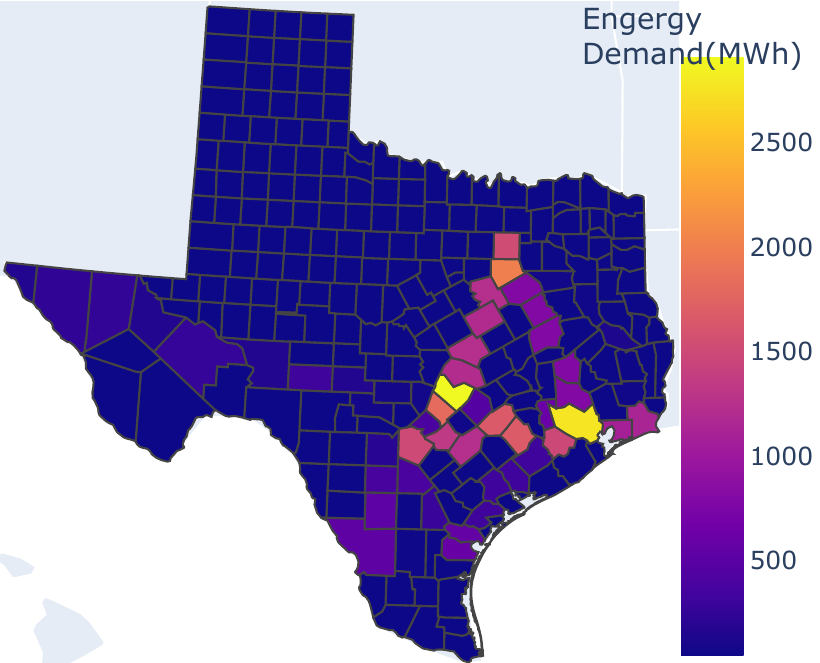}
\vfill
\includegraphics[width=0.75\linewidth]{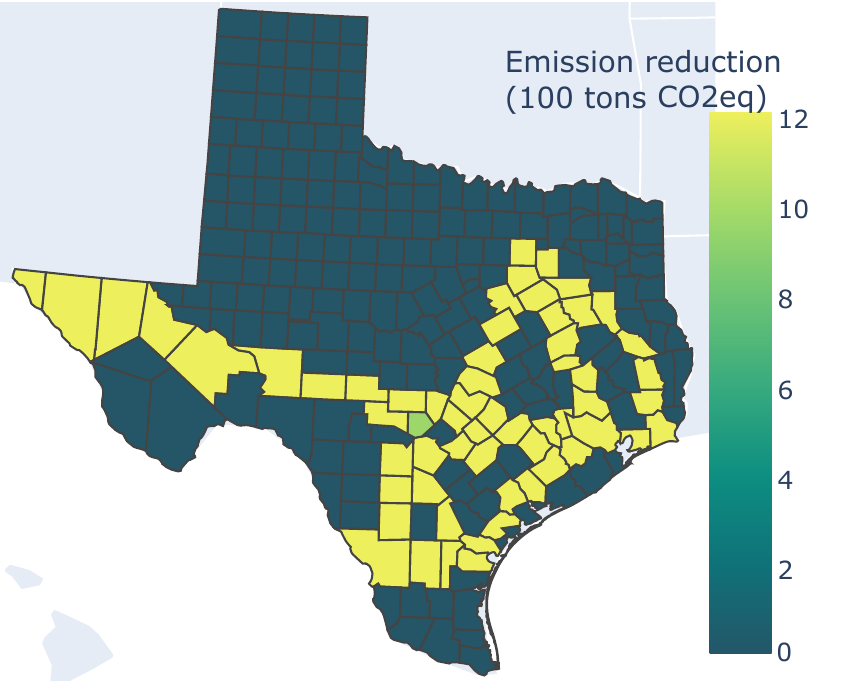}
\caption{Spatial energy demand due to freight electrification (top) and emission reduction due to freight electrification (bottom) per county map for the state of Texas.}
\label{fig:emissionredmap}
\end{figure}

To have a better visualization of distribution of energy demand and emission reduction, for the comparison between BEV fleet and ICEV fleet, we plotted energy demand and the emission reduction due to fleet electrification in Fig.~\ref{fig:emissionredmap}. It is seen that most of the counties the major routes passing through have higher energy demand and emission reduction. Harris County has the highest energy demand followed by Travis County. In contrast, Hays County in the center of the map is showing a lower emission reduction than the neighboring counties, while this county has the largest diesel fuel demand. This is because the BEV fleet also consumes higher electric energy and leads to a higher net-carbon emission.

\begin{figure}[!t]
\centering
\includegraphics[width=0.42\textwidth]{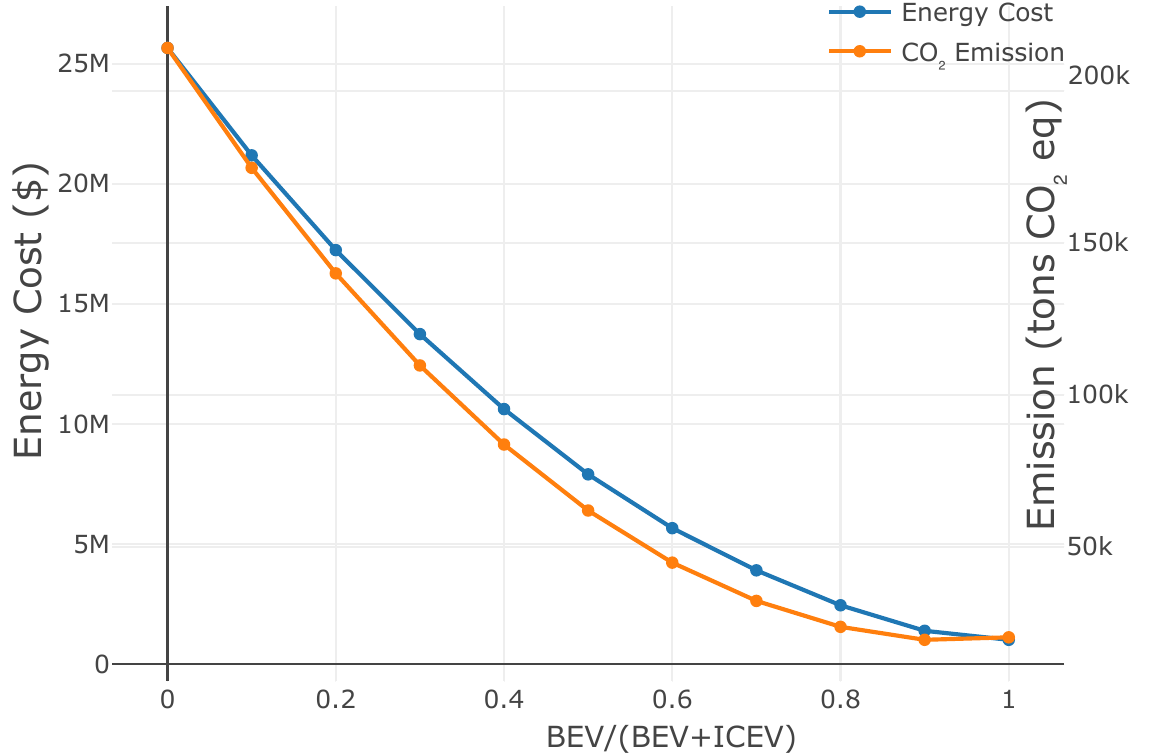}
\caption{Sensitivity analysis of energy cost and $\text{CO}_2$ emission by varying from 0\% to 100\% of BEV in the mixed type of fleet for the case study of Texas.}
\label{fig:bevmix}
\end{figure}

Further analysis of the BEV and ICEV fleet is conducted by mixing the fleet with different types of vehicle and comparing the energy cost and $\text{CO}_2$ emission as we increase the penetration of BEV in the fleet. As can be seen in Fig.~\ref{fig:bevmix}, as the penetration of BEV increase in the mixed type of fleet, both the energy cost and $\text{CO}_2$ emission are decreasing, however, the decreasing rate is also slowing down. We can also observe that $\text{CO}_2$ emission decrease faster than energy cost when the penetration of BEV is lower than 60\%. This can indicate the environmental value of long-haul truck electrification is higher than the economic value when the penetration is lower than 60\%.
Note that the optimization problem solved for the Texas cities includes 28 O-D pairs, hence more than 1600 variables included. The proposed algorithm is able to solve the problem within 360s for the entire year.

\section{Conclusion}\label{sec:conclusion}
In this paper, we developed an optimization framework that integrates data from OpenStreetMaps, utility rates from OPENEI utitlity rate database, and freight volume metrics from FAF5 database, to provide a clearer understanding of the operational costs, spatial distribution of energy demand, and environmental impacts associated with electric and diesel freight trucks.
Given the vast nature and intricacy of the data, we designed our method to be both efficient and capable of analyzing large geographic expanses.


This research holds relevance for stakeholders, such as policymakers and utility companies. Policymakers can derive insights to set effective rates for electric truck charging infrastructure, while utility providers can forecast the regions that might see increased energy demand due to electric truck operations.
To validate and illustrate our methodology, we have conducted a detailed case study focusing on Texas, encompassing major urban centers.

Future work will also consider the up-front cost of freight electrification (including charging infrastructure) and provide a comprehensive comparison between BEV and ICEV fleets through a life-cycle analysis (LCA) method. 


%



\bibliography{references,refs}
\bibliographystyle{IEEEtran}

\end{document}